%
\input mssymb

\hsize=15,4truecm
\vsize=23.5truecm
\mathsurround=1pt

\def\chapter#1{\par\bigbreak \centerline{\bf #1}\medskip}

\def\section#1{\par\bigbreak {\bf #1}\nobreak\enspace}

\def\sqr#1#2{{\vcenter{\hrule height.#2pt
      \hbox{\vrule width.#2pt height#1pt \kern#1pt
         \vrule width.#2pt}
       \hrule height.#2pt}}}

\def\k{\kappa}
\def\o{\omega}

\def\d{\delta}

\def\l{\lambda}
\def\r{\rho}

\def\a{\alpha}
\def\b{\beta}

\def\t{\tau}
\def\g{\gamma}
\def\h{\theta}

\def\n{\eta}

\def\Z{\hbox{\bf Z}}


\def\A{{\cal A}}
\def\B{{\cal B}}

\def\M{{\bf M}}


\def\th #1 #2. #3\par\par{\medbreak{\bf#1 #2.
\enspace}{\sl#3\par}\par\medbreak}
\def\rem #1 #2. #3\par{\medbreak{\bf #1 #2.
\enspace}{#3}\par\medbreak}
\def\proof{{\bf Proof}.\enspace}
\def\sqr#1#2{{\vcenter{\hrule height.#2pt
      \hbox{\vrule width.#2pt height#1pt \kern#1pt
         \vrule width.#2pt}
       \hrule height.#2pt}}}
\def\eop{\mathchoice\sqr34\sqr34\sqr{2.1}3\sqr{1.5}3}

                                                                     %
                                                                     %
\newdimen\refindent\newdimen\plusindent                              %
\newdimen\refskip\newdimen\tempindent                                %
\newdimen\extraindent                                                %
                                                                     %
                                                                     %
\def\ref#1 #2\par{\setbox0=\hbox{#1}\refindent=\wd0                  %
\plusindent=\refskip                                                 %
\extraindent=\refskip                                                %
\advance\extraindent by 30pt                                         %
\advance\plusindent by -\refindent\tempindent=\parindent 
\parindent=0pt\par\hangindent\extraindent 
{#1\hskip\plusindent #2}\parindent=\tempindent}                      %
\refskip=\parindent                                                  %
                                                                     %

\def\ol{\overline}
\def\tr{{\rm tr}}
\def\empty{\emptyset}

\def\raj{\restriction}

\def\ol{\overline}

\def\nda{\mathrel{\lower0pt\hbox to 3pt{\kern3pt$\not$\hss}\downarrow}}
\def\nbot{\mathrel{\lower0pt\hbox to 4pt{\kern3pt$\not$\hss}\bot}}
\def\ekom{\mathrel{\lower3pt\hbox to 0pt{\kern3pt$\sim$\hss}\mapsto}}

\def\anR{\mathrel{\lower1pt\hbox to 2pt{\kern3pt$R$\hss}\not}}
\def\anoR{\mathrel{\lower1pt\hbox to 2pt{\kern3pt$\overline{R}$\hss}\not}}

\def\anRm{\mathrel{\lower1pt\hbox to 2pt{\kern3pt$R^{-}$\hss}\not}}

\def\ndda{\mathrel{\lower0pt\hbox to 1pt{\kern3pt$\not$\hss}\downdownarrows}}

\null
\vskip 2truecm
\centerline{\bf CONSTRUCTING STRONGLY EQUIVALENT NONISOMORPHIC}
\centerline{\bf MODELS FOR UNSUPERSTABLE THEORIES, PART B}
\vskip 1truecm
\centerline{Tapani Hyttinen and Saharon Shelah$^{*}$}
\vskip 2.5truecm

\chapter{Abstract}

We study how equivalent nonisomorphic models of unsuperstable theories
can be.  We measure the equivalence by Ehrenfeucht-Fraisse games.
This paper continues [HS].

\vskip 1.5truecm

\chapter{1. Introduction}

In [HT] we started the studies of so called strong nonstructure theorems.
By strong nonstructure theorem we mean a theorem which says that
if a theory belongs to some class of theories then it has very equivalent
nonisomorphic models. Usually the equivalence is measured by the
length of the Ehrenfeucht-Fraisse games (see Definition 2.2)
in which $\exists$ has
a winning strategy. These theorems are called nonstructure theorems
because intuitively the models must be complicated if they
are very equivalent but still nonisomorphic. Also
structure theorems usually imply that a certain degree of equivalence
gives isomorphism (see f.ex. [Sh1] (Chapter XIII)).

In [HT] we studied mainly unstable theories. We also looked
unsuperstable theories but we were not able
to say much if the equivalence
is measured by the
length of the Ehrenfeucht-Fraisse games in which $\exists$ has
a winning strategy. In this paper we make a new attempt to study
the unsuperstable case.

The main result of this paper is the following:
if $\l =\mu^{+}$,
$cf(\mu )=\mu$, $\k =cf(\k )<\mu$, $\l^{<\k}=\l$, $\mu^{\k}=\mu$
and $T$ is an unsuperstable
theory, $\vert T\vert\le\l$ and $\k (T)>\k$, then there are models
$\A$, $\B\models T$ of cardinality $\l$ such that
$$\A\equiv^{\l}_{\mu\times\k}\B\ \ \hbox{\sl and}
\ \ \A\not\cong\B .$$
In [HS] we proved this theorem in a special case.

\vskip 1.2truecm

\noindent
$^{*}$ Partially supported by the United States Israel Binational
Science Foundation. Publ. 529.

\vfill
\eject

\relax FromTheorem 4.4 in [HS] we get the following theorem easily:
Let $T_{c}$ be the canonical
example of unsuperstable theories i.e. $T_{c}=Th(( ^{\o}\o,E_{i})_{i<\o})$
where $\n E_{i}\xi$ iff for all $j\le i$, $\n (j)=\xi (j)$.

\th 1.1 Theorem. ([HS]) Let $\l =\mu^{+}$ and
$I_{0}$ and $I_{1}$ be models of $T_{c}$ of cardinality $\l$.
Assume $\l\in I[\l ]$.
Then
$$I_{0}\equiv^{\l}_{\mu\times\o +2}I_{1}\ \ \Leftrightarrow
\ \ I_{0}\cong I_{1}.$$

\noindent
So the main result of Chapter 3 is essentially the best possible.

In the introduction of [HT] there is more background for strong
nonstructure theorems.

\chapter{2. Basic definitions}

In this chapter we define the basic concepts we shall use and
construct two linear orders needed in Chapter 3.

\th 2.1 Definition. Let  $\l$ be a cardinal and $\a$ an ordinal.
Let $t$ be a tree (i.e. for all $x\in t$, the set $\{ y\in t\vert\ y<x\}$
is well-ordered by the ordering of $t$).
If $x,y \in t$ and $\{ z \in t \mid z < x \} = \{ z \in t \mid z < y \}$,
then we denote $x \sim y$, and the equivalence class
of $x$ for $\sim$ we denote $[x]$.
By a $\l, \a$-tree $t$ we mean a
tree which satisfies:

(i) $\vert [x] \vert < \l$ for every $x \in t$;

(ii) there are no branches of length $\ge \a$ in $t$;

(iii) $t$ has a unique root;

(iv) if $x,y \in t$, $x$ and $y$ have no immediate predecessors
and $x\sim y$, then
$x=y$.

Note that in a $\l,\a$-tree each ascending sequence of a limit length
has at most one supremum.

\th 2.2 Definition. Let $t$ be a tree and $\k$ a cardinal.
The  Ehrenfeucht-Fraisse game
of length $t$ between models $\A$ and $\B$,
$G^{\k}_{t}(\A, \B)$, is the following.
At each move $\a$:

(i) player $\forall$ chooses $x_\a \in t$, $\k_{\a}<\k$ and
either $a_\a^\b \in \A$, $\b <\k_{\a}$ or $b_\a^\b \in \B$,
$\b <\k_\a$, we will denote this sequence
by $X_{\a}$;

(ii) if $\forall$ chose from $\A$ then
$\exists$ chooses $b_\a^\b \in \B$, $\b <\k_\a$, else
$\exists$ chooses
$a_\a^\b \in \A$, $\b <\k_\a$, we will denote this sequence by $Y_{\a}$.
\medskip
\noindent
$\forall$ must move so that $(x_\b)_{\b \le \a}$
form a strictly increasing sequence in $t$.
$\exists$ must move so that
$\{ (a_\g^\b, b_\g^\b) \vert \g \le \a , \b <\k_\g \}$
is a partial isomorphism from $\A$ to $\B$.
The player who first has to break the rules loses.

We write $\A\equiv^{\k}_{t}\B$ if $\exists$ has a winning strategy
for $G^{\k}_{t}(\A ,\B )$.

\th 2.3 Definition. Let $t$ and $t'$ be trees.

(i) If $x \in t$, then $pred(x)$ denotes the sequence $(x_\a)_{\a < \b}$
of the predecessors of $x$, excluding $x$ itself,
ordered by $<$. Alternatively, we consider $pred(x)$ as a set.
The notation $succ(x)$ denotes the set of immediate
successors of $x$.
If $x,y \in t$ and there is $z$,
such that $x,y \in succ(z)$, then we say that
$x$ and $y$ are brothers.

(ii) By $t^{<\a}$ we mean the set
$$\{ x\in t\vert\ \hbox{\sl the order type of}\ pred(x)
\ \hbox{\sl is}\ <\a\} .$$
Similarly we define $t^{\le\a}$.

(iii) The sum $t \oplus t'$ is defined
as the disjoint union of $t$ and $t'$, except that
the roots are identified.

\th 2.4 Definition. Let $\r_{i}$, $i<\a$,
$\r$ and $\h$ be linear orders.

(i)
We define the ordering $\r\times\h$ as follows:
the domain of
$\r\times\h$ is  $\{ (x,y)\vert\ x\in\r ,\ y\in\h\}$,
and the ordering in  $\r \times \h$ is
defined by last differences, i.e., each point in $\h$ is replaced by a
copy of $\r$;

(ii) We define the ordering $\r +\h$ as follows:
The domain of $\r +\h$ is $(\{ 0\}\times\r )\cup (\{ 1\}\times\h )$
and the ordering in $\r +\h$ is defined by the first difference i.e.
$(i,x)<(j,y)$ iff $i<j$ or $i=j$ and $x<y$.

(iii) We define the ordering $\sum_{i<\a}\r_{i}$ as follows:
The domain of $\sum_{i<\a}\r_{i}$ is $\{ (i,x)\vert\ i\in\a, x\in\r_{i}\}$
and the ordering in $\sum_{i<\a}\r_{i}$ is defined by the first difference i.e.
$(i,x)<(j,y)$ iff $i<j$ or $i=j$ and $x<y$.

\th 2.5 Definition.
We define generalized Ehrenfeucht-Mostowski
models (E-M-models for short).
Let $K$ be a class of models we call index models.
In this definition the notation
$tp_{at}(\ol x,A,\A )$ means the atomic type of $\ol x$
over $A$ in the model $\A$.

Let $\Phi $ be a function. We say that
$\Phi $ is proper for $K$,
if there is a vocabulary $\t_1$ and for
each $I \in K$ a model $\M_1$ and tuples
$\ol a_s$, $s \in I$, of elements of $\M_1$,
such that:

(i) each element in $\M_1$ is an interpretation
of some $\mu(\ol a_{\ol s})$, where $\mu$ is
a $\t_1$-term;

(ii) $tp_{at}(\ol a_{\ol s}, \empty, \M_1)
= \Phi (tp_{at}(\ol s, \empty, I))$.
\medskip
\noindent
Here $\ol s= (s_0,..., s_n)$ denotes a tuple of elements of $I$
and $\ol a_{\ol s}$ denotes
$\ol a_{s_0} \frown \cdots \frown \ol a_{s_n}$.

Note that if $\M_1$, $\ol a_s$, $s \in I$, and
$\M_1'$, $\ol a_s'$, $s \in I$, satisfy
the conditions above, then there is
a canonical isomorphism $\M_1 \cong \M_1'$
which takes $\mu(\ol a_{\ol s})$ in $\M_1$
to $\mu(\ol a'_{\ol s})$ in $\M_1'$.
Therefore we may assume below that
$\M_1$ and $\ol a_s$, $s \in I$, are unique
for each $I$. We denote this unique $\M_1$
by $EM^1(I,\Phi )$ and call it an Ehrenfeucht-Mostowski
model.
The tuples $\ol a_s$, $s \in I$, are
the generating elements of $EM^1(I,\Phi )$,
and the indexed set $(\ol a_s)_{s \in I}$
is the skeleton of $EM^1(I,\Phi )$.

Note that if
$$
tp_{at}(\ol s_1, \empty, I) =
tp_{at}(\ol s_2, \empty, J),
$$
then
$$
tp_{at}(\ol a_{\ol s_1}, \empty, EM^1(I, \Phi )) =
tp_{at}(\ol a_{\ol s_2}, \empty, EM^1(J, \Phi )).
$$

\th 2.6 Definition. Let $\h$ be a linear order and $\k$ infinite regular
cardinal.
Let $K^\k_\tr(\h)$ be the class of models
of the form
$$
I=(M,<,\ll ,H,P_\a )_{\a \le \k},
$$
where $M \subseteq \h^{\le \k}$ and:

(i) $M$ is closed under initial segments;

(ii) $<$ denotes the initial segment relation;

(iii) $H(\n, \nu)$ is the maximal common
initial segment of $\n$ and $\nu$;

(iv) $P_\a = \{ \n \in M \mid length(\n) = \a \}$;

(v)
$\n \ll \nu$ iff either $\n <\nu$ or there is
$n < \k$ such that $\n(n)  < \nu(n)$
and $\n \raj n = \nu \raj n$.
\medskip
\noindent
Let $K^\k_\tr = \bigcup \{ K^\k_\tr(\h) \mid \h$ a linear order $\}$.

If $I \in K^\k_\tr(\h)$ and $\n, \nu \in I$, we define
$\n <_s \nu$ iff $\n$ and $\nu$ are brothers and
$\n < \nu$. But we do not put $<_s$ to the vocabulary of
$I$.

Thus the models in $K^\k_\tr$ are
lexically ordered trees of height $\k + 1$ from which we
have removed the relation $<_s$ and where we have
added relations indicating the levels and a function
giving the maximal common predecessor.

The following theorem gives us means
to construct for $T$ E-M-models such that
the models of $K^\k_\tr$ act as index models.
Furthermore the properties of
the models of $K^\k_\tr$ are reflected to
these E-M-models.

\th 2.7 Theorem. ([Sh1]).
Suppose $\t \subseteq \t_1$,
$T$ is a complete $\t$-theory,
$T_1$ is a complete $\t_1$-theory
with Skolem functions and $T \subseteq T_1$.
Suppose further that $T$ is unsuperstable, $\k(T)>\k$
and $\phi_n(\ol x,\ol y_n)$, $n < \k$,
witness this.
(The definition of witnessing is not needed in this paper.
See [Sh1].)

Then there is a function $\Phi $, which is proper for
$K^\k_\tr$, such that for every $I \in K^\k_\tr$,
$EM^1(I,\Phi )$ is a $\t_1$-model of $T_1$,
for all $\n\in I$, $\ol a_{\n}$ is finite
and for $\n ,\xi \in P_n^I, \nu \in P_\k^I$,

(i) if $I \models \n < \nu$, then
$EM^1(I,\Phi ) \models \phi_n(\ol a_\nu, \ol a_\n)$;

(ii) if $\n$ and $\xi$ are brothers and
$\n <\nu$ then $\xi =\n$ iff $EM^1(I,\Phi ) \models
\phi_n(\ol a_\xi, \ol a_\nu)$.
$\eop$

Above $\phi_n(\ol x, \ol y_n)$
is a first-order $\t$-formula.
We denote the reduct
$$EM^1(I,\Phi ) \raj \t$$
by $EM(I,\Phi )$. In order to simplify the notation, instead of
$\ol a_{\n}$, we just write $\n$. It
will be clear from the context, whether $\n$ means $\ol a_{\n}$
or $\n$.

Next we construct two linear orders needed in the next chapter.
The first of these constructions is a modification of a
linear order construction in [Hu] (Chapter 9).

\th 2.8 Definition. Let $\g$ be an ordinal closed under
ordinal addition and let
$\h_{\g} =(\ ^{<\o}\g ,<)$, where $<$ is defined by $x<y$ iff

(i) $y$ is an initial segment of $x$

\noindent
or

(ii) there is $n<min\{ length(x),length(y)\}$ such that
$x\raj n=y\raj n$ and $x(n)<y(n)$.

\th 2.9 Lemma. Assume $\g$ in an ordinal closed under ordinal addition.
Let $x\in\h_{\g}$, $length(x)=n<\o$ and $\a<\g$.
Let $A^{\a}_{x}$ be the set of all elements $y$ of $\h_{\g}$ which
satisfy:

(i) $x$ is an initial segment of $y$ (not necessarily proper);

(ii) if $length(y)>n$ then $y(n)\ge\a$.

\noindent
Then $(A^{\a}_{x},<\raj A^{\a}_{x})\cong\h_{\g}$.

\proof Follows immediately from the definition of $\h_{\g}$. $\eop$

If $\a\le\b$ are ordinals then by $(\a ,\b ]$ we mean the unique
ordinal order isomor\-phic to
$$\{ \d\vert\ \a <\d\le\b\}\cup\{ \d\vert\ \d =\a\ \hbox{\rm and limit}\}$$
together
with the natural ordering. Notice that if $(\a_{i})_{i<\d}$ is
strictly increasing continuous sequence of ordinals, $\a_{0}=0$,
$\b =sup_{i<\d}\a_{i}$ and for all successor $i<\d$, $\a_{i}$ is
successor,
then $\sum_{i<\d}(\h\times (\a_{i},\a_{i+1}])\cong\h\times\b$,
for all linear-orderings $\h$.

\th 2.10 Lemma. Let $\g$ be an ordinal closed under ordinal
addition and not a cardinal.

(i) Let $\a <\g$ be an ordinal. Then
$$\h_{\g}\cong\h_{\g}\times (\a +1).$$

(ii) Let $\a <\b <\vert\g\vert^{+}$. Then
$$\h_{\g}\cong\h_{\g}\times (\a ,\b ].$$

\proof (i) For all $i<\a$ we let $x_{i}=(i)$. Then by the definition
of $\h_{\g}$,
$$\h_{\g}\cong (\sum_{i<\a}A^{0}_{x_{i}})+A^{\a}_{()},$$
where by $()$ we mean the empty sequence. By Lemma 2.9
$$(\sum_{i<\a}A^{0}_{x_{i}})+A^{\a}_{()}\cong\h_{\g}\times (\a +1).$$

(ii) We prove this by induction on $\b$. For $\b =1$ the claim follows
from (i).
Assume we have proved the claim for $\b <\b '$ and we prove it for
$\b '$. If $\b '=\d +1$, then by induction assumption
$$\h_{\g}\cong\h_{\g}\times (\a ,\d ]$$
and so
$$\h_{\g}\times (\a ,\d +1]\cong\h_{\g}+\h_{\g}\cong\h_{\g}$$
by (i).

If $\b '$ is limit, then we choose a strictly increasing
continuous sequence of
ordinals
$(\b_{i})_{i<cf(\b ')}$, so that $\b_{0}=\a$,
$sup_{i<cf(\b ')}\b_{i}=\b '$ and for all successor $i<cf(\b ')$,
$\b_{i}$ is successor.
Then
$$\h_{\g}\times (\a ,\b ']\cong\sum_{i<cf(\b ')}
(\h_{\g}\times (\b_{i},\b_{i+1}])+\h_{\g}.$$
By induction assumption
$$\sum_{i<cf(\b ')}
(\h_{\g}\times (\b_{i},\b_{i+1}])+\h_{\g}\cong\h_{\g}\times(cf(\b')+1).$$
Because $\g$ is not a cardinal, $cf(\b')<\g$ and so by (i)
$$\h_{\g}\times(cf(\b')+1)\cong\h_{\g}.$$
$\eop$

\th 2.11 Corollary. Let $\g$ be an ordinal closed under ordinal
addition and not a cardinal.
If $\a <\vert\g\vert^{+}$ is
a successor ordinal then $\h_{\g}\cong\h_{\g}\times\a$.

\proof Follows immediately from Lemma 2.10 (ii).
$\eop$

\th 2.12 Lemma. Assume $\mu$ is a regular cardinal and
$\l =\mu^{+}$.
Then there are linear order $\h$ of power $\l$,
one-one and onto function
$h:\h\rightarrow\l\times\h$
and order isomorphisms $g_{\a}:\h\rightarrow\h$ for $\a <\l$
such that the
following holds:

(i) if $g_{\a}(x)=y$ then $x\ne y$ and either

(a) $h(x)=(\a ,y)$

or

(b) $h(y)=(\a ,x)$

but not both,

(ii) if for some $x\in\h$, $g_{\a}(x)=g_{\a '}(x)$
then $\a =\a '$,

(iii) if $h(x)=(\a ,y)$ then $g_{\a}(x)=y$ or $g_{\a}(y)=x$.

\proof Let the universe of $\h$ be $\mu\times\l$. The ordering
will be defined by induction.
Let
$$f:\l\rightarrow\l\times\l$$
be one-one, onto and if $\a <\a'$, $f(\a )=(\b ,\g )$ and $f(\a')=(\b',\g')$
then $\g <\g'$.
This $f$ is used only to guarantee that
in the induction we pay
attention to every $\b <\l$ cofinally often.

By induction on $\a <\l$ we do the following:
Let $f(\a )=(\b ,\g )$. We define $\h^{\a}=
(\mu\times (\a+1),<^{\a})$,
$h^{\a}:\h^{\a}\rightarrow \l\times\h^{\a}$
and order isomorphisms (in the ordering $<^{\a}$)
$$g^{\a}_{\b}:\h^{\a}\rightarrow\h^{\a}$$
so that

(i) if $\a <\a'$ then $h^{\a}\subseteq h^{\a '}$ and
$<^{\a}\subseteq <^{\a '}$,

(ii) if $\a <\a'$, $f(\a )=(\b ,\g )$ and $f(\a ')=(\b ,\g ' )$
then $g^{\a}_{\b}\subseteq g^{\a '}_{\b}$,

(iii) if $g^{\a}_{\b}(x)=y$ then $x\ne y$ and either

(a) $h^{\a}(x)=(\b ,y)$

or

(b) $h^{\a}(y)=(\b ,x)$

but not both.

The induction is easy since at each stage we have $\mu$
"new" elements to use:
Let $B\subseteq\mu\times\a$ be the set of those element
from $\mu\times\a$ which are not in the domain of any $g^{\a '}_{\b}$
such that $\a '<\a$ and $f(\a ')=(\b ,\g ')$ for some $\g '$.
(Notice that $B$ is also the set of those element
from $\mu\times\a$ which are not in the range of any $g^{\a '}_{\b}$
such that $\a '<\a$ and $f(\a ')=(\b ,\g ')$ for some $\g '$.)
Clearly if $B\ne\empty$ then $\vert B\vert =\mu$.

Let $A_{i}$, $i\in\Z$, be a partition
of $\mu\times\{ \a\}$ into sets of power $\mu$.
We first define $g^{\a}_{\b}$ so that the following is true:

(a) $g^{\a}_{\b}$ is one-one,

(b) if $B\ne\empty$ then $g^{\a}_{\b}\raj A_{0}$ is onto $B$
otherwise $g^{\a}_{\b}\raj A_{0}$ is onto $A_{-1}$,

(c) if $B\ne\empty$ then $g^{\a}_{\b}\raj B$ is onto $A_{-1}$,

(d) for all $i\ne 0$, $g^{\a}_{\b}\raj A_{i}$ is onto $A_{i-1}$.

By an easy induction on $\vert i\vert<\o$ we can
define $<^{\a}$ so that $<^{\a'}\subseteq <^{\a}$ for all
$\a '<\a$ and $g^{\a}_{\b}$ is an order isomorphism.
We define the function
$h^{\a}\raj (\mu\times\{\a\} )$ as follows:

(a) if $B=\empty$ then $h^{\a}(x)=(\b ,g^{\a}_{\b}(x))$,

(b) if $B\ne\empty$ and $i\ge 0$ and $x\in A_{i}$ then
$h^{\a}(x)=(\b ,g^{\a}_{\b}(x))$,

(c) if $B\ne\empty$ and $i<0$ and $x\in A_{i}$ then
$h^{\a}(x)=(\b ,y)$ where
$y\in A_{i+1}$ or $B$ is the unique element such that
$g^{\a}_{\b}(y)=x$.

\noindent
It is easy to see that (iii) above is satisfied.

We define $\h=(\mu\times\l ,<)$, where
$<=\bigcup_{\a <\l}<^{\a}$, $h=\bigcup_{\a <\l}h^{\a}$
and for all $\b <\l$
we let
$g_{\b}=\bigcup\{ g^{\a}_{\b}\vert\ \a <\l ,\ f(\a )=(\b ,\g )
\ \hbox{\rm for some}\ \g\}$.
Clearly these satisfy (i). (ii) follows from the fact that
if $g^{\a}_{\b}(x)=y$ then either $x\in\mu\times\{\a\}$ and
$y\in\mu\times(\a +1)$
or $y\in\mu\times\{\a\}$ and $x\in\mu\times(\a +1)$.
(iii) follows immediately from the definition of $h$.
$\eop$

\chapter{3. On nonstructure of unsuperstable theories}

In this chapter we will prove the main theorem of this paper i.e.
Conclusion 3.19. The idea of the proof continues
III Claim 7.8 in [Sh2].
Throughout this chapter we assume that
$T$ is an unsuperstable theory, $\vert T\vert <\l$ and $\k (T)>\k$.
The
cardinal assumptions are: $\l =\mu^{+}$,
$cf(\mu )=\mu$, $\k =cf(\k )<\mu$, $\l^{<\k}=\l$, $\mu^{\k}=\mu$.

If $i<\k$ we say that $i$ is
of type $n$, $n=0,1,2$, if there are a limit ordinal $\a <\k$ and
$k<\o$ such that $i=\a +3k+n$. 

We define linear orderings $\h_{n}$, $n<3$, as follows.
Let $\h_{0}=\l$ and $\h_{1}$, $h'$ and $g_{\a}$, $\a <\l$, as
$\h$, $h$ and $g_{\a}$ in Lemma 2.12.
Let $\h_{2}=\h_{\mu\times\o}\times\l$, where $\h_{\mu\times\o}$
is as in Definition 2.8.

For $n<2$, let $J^{-}_{n}$ be the set of sequences $\n$ of length $<\k$
such that

(i) $\n\ne ()$;

(ii) $\n(0)=n$;

(iii) if $0<i<length(\n )$ is of type $m<3$ then
$\n (i)\in\h_{m}$.

Let
$$f:(\l -\{ 0\} )\rightarrow\{ (\n ,\xi )\in J^{-}_{0}\times J^{-}_{1}\vert
\ length(\n )=length(\xi )\ \hbox{\rm is of type}\ 1\}$$
be one-one and onto.
Then we define
$$h:\h_{1}\rightarrow J^{-}_{0}\cup J^{-}_{1}$$
and order isomorphisms
$$g_{\n ,\xi}:succ(\n )\rightarrow succ(\xi ),$$
for $(\n ,\xi )\in rng(f)$, as follows:

(i) $g_{\n ,\xi}(\n\frown (x))=\xi\frown (g_{\a}(x))$, where $\a$ is
the unique ordinal such that
$f(\a )=(\n ,\xi )$;

(ii) Assume $h'(x)=(\a ,y)$, $\a\ne 0$, and $f(\a )=(\n ,\xi )$. Then
$h(x)=\xi\frown (y)$ if $g_{\a}(x)=y$ otherwise
$h(x)=\n\frown (y)$. If $h'(x)=(0,y)$ then $h(x)=(0)$ (here the
idea is to define $h(x)$ so that $length(h(x))$ is not of type 2).

\th 3.1 Lemma.
Assume $\n\in J^{-}_{0}$ and $\xi\in J^{-}_{1}$ are such that
$m=length(\n )=length(\xi )$ is of type $2$. Let $m=n+1$.
If $g_{\n ,\xi}(\n ')=\xi '$ then either

(a) $h(\n '(n))=\xi '$

or

(b) $h(\xi '(n))=\n '$

but not both.

\proof We show first that either (a) or (b)
holds. So we
assume that (a) is not true and prove that (b) holds.
Let $\n '(n)=x$, $\xi '(n)=y$
and $f(\a )=(\n ,\xi )$. Now $g_{\a}(x)=y$, $x\ne y$ and either
$h'(x)=(\a ,y)$ or $h'(y)=(\a ,x)$. Because (a) is not true
$h'(x)\ne (\a ,y)$ and so $h'(y)=(\a ,x)$.
We have two cases:

(i) Case $y>x$: Because $g_{\a}$ is order-precerving,
$g_{\a}(y)>y>x$. So $g_{\a}(y)\ne x$ and by
the definition of $h$, $h(y)=\n\frown (x)=\n '$.

(ii) Case $y<x$: As the case $y>x$.

Next we show that it is impossible that both (a) and (b) holds.
For a contradiction assume that this is not the case. Then
(a) implies that there is $\b$ such that $h'(x)=(\b ,y)$
and $g_{\b}(x)=y$. On the other hand (b) implies that there
is $\g$ such that $h'(y)=(\g ,x)$ and $g_{\g}(y)\ne x$.
By Lemma 2.12 (iii), $g_{\g}(x)=y$.
By Lemma 2.12 (ii) $\b =\g$.
So $h'(y)=(\b ,x)$ and $h'(x)=(\b ,y)$, which contradicts Lemma 2.12 (i).
$\eop$

For $n<2$, let $J^{+}_{n}$ be the set of sequences $\n$ of length $\le\k$
such that

(i) $\n\ne ()$;

(ii) $\n(0)=n$;

(iii) if $0<i<length(\n )$ is of type $m<3$ then
$\n (i)\in\h_{m}$.

Let $e:\h_{1}\rightarrow\l$ be one-one and onto.
We define functions $s$ and $d$
as follows: if $i<length(\n )$ is of type $0$
then $d(\n ,i)=\n (i)$ and
$s(\n ,i)=\n (i)$,
if $i<length(\n )$ is of type $1$
then $d(\n ,i)=\n (i)$ and
$s(\n ,i)=e(\n (i))$ and
if $i<length(\n )$ is of type $2$ and
$\n (i)=(d,s)$ then $d(\n ,i)=d$ and $s(\n ,i)=s$.

For $n<2$ and $\g<\l$, we define
$$J^{+}_{n}(\g )=\{ \n\in J^{+}_{n}\vert\ \hbox{\rm for all}
\ i<length(\n ),\ s(\n ,i)<\g\} ,$$
$J^{-}_{n}(\g )=J^{+}_{n}(\g )\cap J^{-}_{n}$.

Let us fix $d\in\h_{1}$ so that
$h(d)=(0)$.

\th 3.2 Definition. For all $\n\in J^{-}_{0}$ and
$\xi\in J^{-}_{1}$ such that $n=length(\n )=length(\xi )$ is
of type 1,
let $\a (\n ,\xi )$ be the set of ordinals $\a <\l$ such that
for all $\n '\in succ(\n )$, $s(\n ',n)<\a$ iff
$s(g_{\n ,\xi}(\n '),n)<\a$ and $e(d)<\a$. Notice that $\a (\n ,\xi )$
is a closed and unbounded subset of $\l$.
By $\a (\b )$, $\b <\l$, we mean
$$Min\ \bigcap\{ \a (\n ,\xi)\vert\ \n\in J^{-}_{0}(\b ),
\ \xi\in J^{-}_{1}(\b ),\ length(\n )=length(\xi )
\ \hbox{\sl is of type}\ 1\} .$$

\th 3.3 Definition. For all $\n\in J^{+}_{0}$ and $\xi\in J^{+}_{1}$,
we write $\n R^{-}\xi$ and $\xi R^{-}\n$ iff

(i) $\n (j)=\xi (j)$ for all
$0<j<min\{ length(\n ),length(\xi )\}$ of type 0;

(ii) for all $j<min\{ length(\n ),length(\xi )\}$ of type 1
$\xi\raj (j+1)=g_{\n\raj j ,\xi\raj j}(\n\raj (j+1))$.

Let $length(\n )=length(\xi )=j+1$, $j$ of type 1,
and $\n R^{-}\xi$.
We write $\n\rightarrow\xi$ if
$h(\n (j))=\xi$.
We write $\xi\rightarrow\n$ if
$h(\xi (j))=\n$.

\th 3.4 Remark. If $\xi\rightarrow\n$
and $\xi\rightarrow\n '$ then $\n =\n '$ and if $\n R^{-}\xi$ then
$\n\rightarrow\xi$ or $\xi\rightarrow\n$ but not both.

\th 3.5 Definition. Let $\n\in J_{0}^{+}-J_{0}^{-}$ and
$\xi\in J_{1}^{+}-J_{1}^{-}$. We write $\n R\xi$ and $\xi R\n$ iff

(i) $\n R^{-}\xi$;

(ii) for every $j<\k$ of type 2, $\n$ and $\xi$ satisfy the following:
if $\n\raj j\rightarrow\xi\raj j$ then $s(\n ,j)\le s(\xi ,j)$ and
if $\xi\raj j\rightarrow\n\raj j$ then $s(\xi ,j)\le s(\n ,j)$;

(iii) the set $W^{\k}_{\n ,\xi}$ is bounded in $\k$, where
$W^{\k}_{\n ,\xi}$ is defined in the following way:
Let $\n\in J^{+}_{0}-J^{<\d}_{0}$ (see Definition 2.3 (ii))
and $\xi\in J^{+}_{1}-J^{<\d}_{1}$ then
$$W^{\d}_{\n ,\xi}=W^{\d}_{\xi ,\n}=V^{\d}_{\n ,\xi}\cup U^{\d}_{\n ,\xi},$$
where
$$V^{\d}_{\n ,\xi}=\{ j<\d\vert\ j\ \hbox{\sl is of type 2 and}
\ \xi\raj j\rightarrow\n\raj j\ \hbox{\sl and}$$
$$\ cf(s(\n ,j))=
\mu\ \hbox{\sl and}
\ s(\xi ,j)=s(\n ,j)\}$$
and
$$U^{\d}_{\n ,\xi}=\{ j<\d\vert\ j\ \hbox{\sl is of type 2 and}
\ \n\raj j\rightarrow\xi\raj j\ \hbox{\sl and}$$
$$cf(s(\xi ,j))=
\mu\ \hbox{\sl and}
\ s(\n ,j)=s(\xi ,j)\} .$$

Our next goal is to prove that if $J_{0}$ and $J_{1}$ are such that

(i) $J^{-}_{n}\subseteq J_{n}\subseteq J^{+}_{n}$, $n=0,1$
\noindent
and

(ii) if $\n\in J^{+}_{0}$, $\xi\in J^{+}_{1}$ and $\n R\xi$ then
$\n\in J_{0}$ iff $\xi\in J_{1}$,

\noindent
then $(J_{0}, <,<_{s})\equiv^{\l}_{\mu\times\k}
(J_{1},<,<_{s})$, where $<$ is the initial segment relation
and $<_{s}$ is the union of natural orderings of $succ(\n )$ for all
elements $\n$ of the model.
\relax Fromnow on
in this chapter we assume
that $J_{0}$ and $J_{1}$ satisfy (i) and (ii) above.

The relation $R$ designed not only to guarantee the equivalence
but also to make it possible to prove that the final models
are not isomorphic. Here (iii) in the definition of $R$
plays a vital role. The pressing down elements $\n$ such that
$cf(s(\n ,i))=\mu$, $i$ of type 2,
in (iii) prevents
us from adding too many elements to $J_{n}-J^{-}_{n}$, $n<2$.

For $n<2$, we write
$J_{n}(\g )=J^{+}_{n}(\g )\cap J_{n}$.

\th 3.6 Definition. Let $\a <\k$. $G_{\a}$ is the family of all
partial functions $f$ satisfying:

(a) $f$ is a partial isomorphism from $J_{0}$ to $J_{1}$;

(b) $dom(f)$ and $rng(f)$ are closed under initial segments and
for some $\b <\l$ they are
included in $J_{0}(\b )$ and
$J_{1}(\b )$, respectively;

(c) if $f(\n )=\xi$ then $\n R^{-}\xi$;

(d) if $\n\in J^{+}_{0}$, $\xi\in J^{+}_{1}$, $f(\n )=\xi$ and
$j<length(\n )$ of type 2, then $\n$ and $\xi$ satisfy the following:
if $\n\raj j\rightarrow\xi\raj j$ then $s(\n ,j)\le s(\xi ,j)$ and
if $\xi\raj j\rightarrow\n\raj j$ then $s(\xi ,j)\le s(\n ,j)$;

(e) assume $\n\in J^{+}_{0}-J^{<\d}_{0}$ and
$\{\n\raj\g\vert\ \g <\d\}\subseteq dom(f)$ and
let
$$\xi =\bigcup_{\g <\d}f(\n\raj\g ),$$
then
$W^{\d}_{\n ,\xi}$
has order type $\le\a$;

(f) if $\n\in dom(f)$ and $length(\n )$ is of type $2$ then
$$\{ i<\l\vert\ \hbox{\sl for all}\ d\in\h_{2} ,
\ \n\frown ((d,i))\in dom(f)\} =$$
$$\{ i<\l\vert\ \hbox{\sl for some}\ d\in\h_{2} ,
\ \n\frown ((d,i))\in dom(f)\} =$$
$$\{ i<\l\vert\ \hbox{\sl for all}\ d\in\h_{2} ,
\ f(\n )\frown ((d,i))\in rng(f)\} =$$
$$\{ i<\l\vert\ \hbox{\sl for some}\ d\in\h_{2} ,
\ f(\n )\frown ((d,i))\in rng(f)\}$$
is an ordinal.

We define $F_{\a}\subseteq G_{\a}$ by replacing (f) above by

(f') if $\n\in dom(f)$ and $length(\n )$ is of type $2$ then
$$\{ i<\l\vert\ \hbox{\sl for all}\ d\in\h_{2} ,
\ \n\frown ((d,i))\in dom(f)\} =$$
$$\{ i<\l\vert\ \hbox{\sl for some}\ d\in\h_{2} ,
\ \n\frown ((d,i))\in dom(f)\} =$$
$$\{ i<\l\vert\ \hbox{\sl for all}\ d\in\h_{2} ,
\ f(\n )\frown ((d,i))\in rng(f)\} =$$
$$\{ i<\l\vert\ \hbox{\sl for some}\ d\in\h_{2} ,
\ f(\n )\frown ((d,i))\in rng(f)\}$$
is an ordinal
and of cofinality $<\mu$.

The idea in the definition above is roughly the following:
If $f\in G_{\a}$ and $f(\n )=\xi$ then $\n R\xi$ and
the order type of $W^{\d}_{\n ,\xi}$ is $\le\a$.
If $f\in F_{\a}$ then not only $f\in G_{\a}$ but
$f$ is such that for all small $A\subset J_{0}\cup J_{1}$
we can find $g\supset f$ such that $A\subset dom(g)\cup rng(g)$
and $g\in F_{\a}$.

\th 3.7 Definition. For $f,g\in G_{\a}$ we write $f\le g$ if
$f\subseteq g$ and if $\g <\d \le\k$,
$\n\in J^{+}_{0}-J^{<\d}_{0}$, $\n\raj\g\in dom(f)$,
$\n\raj (\g+1)\not\in dom(f)$,
$\n\raj j\in dom(g)$ for all $j<\d$
and $\xi =\bigcup_{j<\d}g(\n\raj j)$,
then $W^{\g}_{\n ,\xi}=W^{\d}_{\n ,\xi}$.

Notice that $f\le g$ is a transitive relation.

\th 3.8 Remark. Let $f\in G_{\a}$. We define
$\overline{f}\supseteq f$ by
$$dom(\overline{f})=dom(f)\cup\{\n\in J_{0}\vert\ \n\raj\g\in
dom(f)\ \hbox{\sl for all}\ \g <length(\n )$$
$$\hbox{\sl and}
\ length(\n )\ \hbox{\sl is limit}\}$$
and if
$\n\in dom(\overline{f})-dom(f)$ then
$$\overline{f}(\n )=\bigcup_{\g <length(\n )}f(\n\raj\g).$$
If $f\in F_{\a}$ then $\overline{f}\in F_{\a}$ and
if $f\in G_{\a}$ then $\overline{f}\in G_{\a}$.

\th 3.9 Lemma. Assume $\a <\k$, $\d\le\mu$, $f_{i}\in F_{\a}$ for all
$i<\d$ and $f_{i}\le f_{j}$ for all $i<j<\d$.

(i) $\bigcup_{i<\d}f_{i}\in G_{\a}$.

(ii) If $\d <\mu$ then
$\bigcup_{i<\d}f_{i}\in F_{\a}$ and $f_{j}\le\bigcup_{i<\d}f_{i}$
for all $j\le\d$.

\proof (i) We have to check that $f=\bigcup_{i<\d}f_{i}$
satisfies (a)-(f) in Definition 3.6. Excluding purhapse (e),
all of these are trivial.

Without loss of generality we may assume
$\d$ is a limit ordinal.
So assume $\n\in J^{+}_{0}-J^{<\b}_{0}$ and
$\{\n\raj\g\vert\ \g <\b\}\subseteq dom(f)$ and
let
$$\xi =\bigcup_{\g <\b}f(\n\raj\g ).$$
We need to show that $W^{\b}_{\n ,\xi}\le\a$.

If there is $i<\d$ such that $\n\raj\g\in dom(f_{i})$ for
all $\g <\b$ then the claim follows immediately from
the assumption $f_{i}\in F_{\a}$. Otherwise
for all $\g <\b$ we let $i_{\g}<\d$ be the least ordinal such that
$\n\raj\g\in dom(f_{i_{\g}})$.
Let $\g^{*}<\b$ be the least ordinal such that $i_{\g^{*}+1}>i_{\g^{*}}$.
Because for all $\g <\b$, $f_{i_{\g}}\in F_{\a}$, we get
$W^{\g}_{\n\raj\g ,\xi\raj\g}$ has order type $\le\a$.
If $\g^{*}<\g'<\b$ then $f_{i_{\g^{*}}}\le f_{i_{\g '}}$ and so
$W^{\g^{*}}_{\n\raj\g^{*} ,\xi\raj\g^{*}}=
W^{\g'}_{\n\raj\g' ,\xi\raj\g'}$.
Because $W^{\b}_{\n ,\xi}=\bigcup_{\g <\b}W^{\g}_{\n\raj\g ,\xi\raj\g}$,
we get $W^{\b}_{\n ,\xi}\le\a$.

(ii) As (i), just check the definitions.
$\eop$

\th 3.10 Lemma. If $\d<\k$, $f_{i}\in G_{i}$ for all $i<\d$ and
$f_{i}\subseteq f_{j}$ for all $i<j<\d$ then
$$\bigcup_{i<\d}f_{i}\in G_{\d}.$$

\proof Follows immediately from the definitions. $\eop$

\th 3.11 Lemma. If $f\in F_{\a}$ and $A\subseteq J_{0}\cup J_{1}$,
$\vert A\vert <\l$, then there is $g\in F_{\a}$ such that $f\le g$
and $A\subseteq dom(g)\cup rng(g)$.

\proof We may assume that $A$ is closed under initial segments.
Let $A'=A\cap (J^{-}_{0}\cup J^{-}_{1})$.
We enumerate $A'=\{ a_{i}\vert\ 0<i<\mu\}$ so that if $a_{i}$ is an
initial segment of $a_{j}$ then $i<j$. Let $\g <\l$ be such that
$A\cup dom(f)\cup rng(f)\subseteq J_{0}(\g )\cup J_{1}(\g )$.
By induction on $i<\mu$ we define functions $g_{i}$.

If $i=0$ we define $g_{i}=f\cup\{ ((0),(1))\}$.

If $i<\mu$ is limit then we define
$$g_{i}=\overline{\bigcup_{j<i}g_{j}}.$$

If $i=j+1$ then there are two different cases.
For simplicity
we assume $a_{i}\in J_{0}$.

(i) $n=length(a_{i})$ is of type $0$ or $1$: Then we choose $g_{i}$ to be such
that

(a) $g_{j}\le g_{i}$;

(b) $g_{i}\in F_{\a}$;

(c) if $\xi\in dom(g_{i})-dom(g_{j})$ then $\xi\in succ(a_{i})$;

(d) if $\xi\in succ(a_{i})$ and $s(\xi ,n)<\g$ then $\xi\in dom(g_{i})$;

(e) if $\xi\in succ(g_{j}(a_{i}))$ and $s(\xi ,n)<\g$ then $\xi\in rng(g_{i})$.

\noindent
Trivially such $g_{i}$ exists.

(ii) $n=length(a_{j})$ is of type $2$: Then we choose $g_{i}$ to be such
that (a)-(c) above and (d')-(f') below are satisfied.

Let
$$\b=sup\{ i+1<\l\vert\ \hbox{\sl for all}\ d\in\h_{2} ,
\ a_{i}\frown ((d,i))\in dom(g_{j})\} .$$

(d') if $\xi\in succ(a_{i})$ then $s(\xi ,n)<\g +2$ iff $\xi\in dom(g_{i})$;

(e') if $\xi\in succ(g_{j}(a_{i}))$ then
$s(\xi ,n)<\g +2$ iff $\xi\in rng(g_{i})$;

(f') $g_{i}\raj\{ \n\in succ(a_{i})\vert\ \b\le s(\n ,n)<\g+1\}$ is
an order isomorphism to
$\{ \n\in succ(g_{j}(a_{i}))\vert\ \b\le s(\n ,n)<\b+1\}$ and
$g_{i}\raj\{ \n\in succ(a_{i})\vert\ \g +1\le s(\n ,n)<\g+2\}$ is
an order isomorphism to
$\{ \n\in succ(g_{j}(a_{i}))\vert\ \b +1\le s(\n ,n)<\g+2\}$.

\noindent
By Corollary 2.11 it is easy to satisfy (d')-(f').
Because $g_{j}\in F_{\a}$, $cf(\b )<\mu$ and
we do not have problems with (a) and (b). So there is $g_{i}$
satisfying (a)-(c) and (d')-(f').

Finally we define
$$g=\overline{\bigcup_{i<\mu}g_{i}}.$$
It is easy to see that $g$ is as wanted (notice that
$f\le g$ follows from the construction, not from Lemma 3.9).
$\eop$

\th 3.12 Lemma. If $f\in G_{\a}$ and $A\subseteq J_{0}\cup J_{1}$,
$\vert A\vert <\l$, then there is $g\in F_{\a +1}$ such that
$f\subseteq g$ and
$A\subseteq dom(g)\cup rng(g)$.

\proof Essentially as the proof of Lemma 3.11. $\eop$

\th 3.13 Theorem. If $J_{0}$ and $J_{1}$ are such that

(i) $J^{-}_{n}\subseteq J_{n}\subseteq J^{+}_{n}$, $n=0,1$
\noindent
and

(ii) if $\n R\xi$, $\n\in J^{+}_{0}$ and $\xi\in J^{+}_{1}$ then
$\n\in J_{0}$ iff $\xi\in J_{1}$,

\noindent
then $(J_{0},<,<_{s})\equiv^{\l}_{\mu\times\k}(J_{1},<,<_{s})$.

\proof Because $\empty\in F_{0}$,
the theorem follows from the previous lemmas. $\eop$

\th 3.14 Corollary. If $J_{0}$ and $J_{1}$ are as above and $\Phi$
is proper for $T$, then
$$EM(J_{0},\Phi )\equiv^{\l}_{\mu\times\k}EM(J_{1},\Phi ).$$

\proof Follows immediately from the definition of
E-M-models and Theorem 3.13. $\eop$

In the rest of this chapter we show that there are trees $J_{0}$ and
$J_{1}$ which satisfy the assumptions of Corollary 3.14 and
$$EM(J_{0},\Phi )\not\cong EM(J_{1},\Phi ).$$

\th 3.15 Lemma. (Claim 7.8B [Sh2]) There are closed increasing
cofinal sequences $(\a_{i})_{i<\k}$
in $\a$,
$\a <\l$ and $cf(\a )=\k$, such that if $i$ is successor then
$cf(\a_{i})=\mu$ and
for all cub $A\subseteq\l$ the set
$$\{\a <\l\vert\ cf(\a )=\k\ \hbox{\sl and}
\ \{\a_{i}\vert\ i<\k\}\subseteq A\cap\a\ \}$$
is stationary.

We define $J_{0}-J^{-}_{0}$ and $J_{1}-J^{-}_{1}$ by using Lemma
3.15. For all $\a <\l$ we define $I^{\a}_{0}$ and $I^{\a}_{1}$.
Let $I^{0}_{0}=J^{-}_{0}$ and $I^{0}_{1}=J^{-}_{1}$.
If $0<\a <\l$, $cf(\a )=\k$,
and there are sequence $(\b_{i})_{i<\k}$ and
$\n\in J^{+}_{0}-J^{-}_{0}$ such that

(i) $(\b_{i})_{i<\k}$ is properly increasing and cofinal in $\a$;

(ii) for all $i<\k$, $cf(\b_{i+1})=\mu$,
$\b_{i+1}>\a (\b_{i})$
and
$\b_{i}\in\{\a_{i}\vert\ i<\k\}$;

(iii) for all $0<i<\k$ of type $0$ or $2$, $s(\n ,i)=\b_{i}$;

(iv) for all $i<\k$ of type 1, $\n (i)=d$;

\noindent
then we choose some such $\n$, let it be $\n_{\a}$, and
define $I^{\a}_{0}$ and $I^{\a}_{1}$ to be the least
sets such that

(i) $\{\n_{\a}\}\cup\bigcup_{\b <\a}I^{\b}_{0}\subseteq I^{\a}_{0}$ and
$\bigcup_{\b <\a}I^{\b}_{1}\subseteq I^{\a}_{1}$

(ii) $I^{\a}_{0}\cup I^{\a}_{1}$ is closed under $R$.

\noindent
Otherwise we let $I^{\a}_{0}=\bigcup_{\b <\a}I^{\b}_{0}$ and
$I^{\a}_{1}=\bigcup_{\b <\a}I^{\b}_{1}$. Finally we define
$J_{0}=\bigcup_{\a <\l}I^{\a}_{0}$ and $J_{1}=\bigcup_{\a <\l}I^{\a}_{1}$.

\th 3.16 Lemma. For all $\a <\l$ and
$\n\in (J_{0}\cup J_{1})-(J^{-}_{0}\cup J^{-}_{1})$, the following
are equivalent:

(i) $\n\in (I^{\a}_{0}\cup I^{\a}_{1})-(\bigcup_{\b <\a}I^{\b}_{0}
\cup\bigcup_{\b <\a}I^{\b}_{1})$.

(ii) $sup\{ s(\n ,i)\vert\ i<\k\} =\a$.

\proof By the construction it is enough to show that (i) implies
(ii). So assume (i). Because of levels of type 0, it is enough
to show that for all $i<\k$,
$s(\n ,i)<\b_{i+1}$. We prove this by induction on $i<\k$.
If $i$ is of type 0, the claim is clear.
If $i$ is of type 1 this follows from
$\b_{i+1}>\a (\b_{i})$ and $e(d)<\a (\b_{i})$ together with the
induction assumption.
For $i$ is of type 2, $i=j+1$,
it is enough to show
that
$s(\n_{\a},i)\ge s(\n ,i)$.
This follows easily from the fact that $\n_{\a}(j)=d$ and
$length(h(d))\ne i$.
$\eop$

\th 3.17 Definition. Let $g:EM(J_{0},\Phi )\rightarrow EM(J_{1},\Phi )$
be an isomorphism. We say that $\a <\l$ is $g$-saturated iff for
all $\n\in J_{0}$ and $\xi_{0},...,\xi_{n}\in J_{1}$ the following
holds: if

(i) $length(\n )=l+1$ and for all $i<l$, $s(\n ,i)<\a$;

(ii) for all $k\le n$ and $i<length(\xi_{k})$, $s(\xi_{k},i)<\a$;

(iii) $g(\n )=t(\d_{0},...,\d_{m})$,
for some term $t$ and $\d_{0},...,\d_{m}\in J_{1}$;

\noindent
then there are $\n '\in J_{0}$ and $\d '_{0},...,\d '_{n}\in J_{1}$
such that

(a) $g(\n ')=t(\d '_{0},...,\d '_{m})$;

(b) $length(\n ')=l+1$ and $\n '\raj l=\n\raj l$;

(c) $s(\n ',l)<\a$;

(d) the basic type of $(\xi_{0},...,\xi_{n},\d_{0},...,\d_{m})$ in
$(J_{1},<,\ll ,H,P_{j})$ is the same as the basic type of
$(\xi_{0},...,\xi_{n},\d '_{0},...,\d '_{m})$.

Notice that for all isomorphisms
$g:EM(J_{0},\Phi )\rightarrow EM(J_{1},\Phi )$ the set of
$g$-saturated ordinals is unbounded in $\l$
and closed under increasing
sequences of length $\a <\l$ if $cf(\a )> \k$.

\th 3.18 Lemma. Let $\Phi$ be proper for $T$. Then
$$EM(J_{0},\Phi )\not\cong EM(J_{1},\Phi ).$$

\proof We write $\A_{\g}$ for the submodel of $EM(J_{0},\Phi )$
generated (in the extended language) by $J_{0}(\g )$. Similarly,
we write $\B_{\g}$ for the submodel of $EM(J_{1},\Phi )$ generated
by $J_{1}(\g )$.
Let $g$ be an one-one function from
$EM(J_{0},\Phi )$ onto $EM(J_{1},\Phi )$.
We say that $g$ is closed in $\g$, if $\A_{\g}\cup\B_{\g}$ is
closed under $g$ and $g^{-1}$.

For a contradiction we assume that $g$ is an isomorphism from
$EM(J_{0},\Phi )$ to $EM(J_{1},\Phi )$.
By Lemma 3.15 we choose $\a <\l$ to be
such that

(i) $cf(\a )=\k$, for all $i<\k$, $g$ is closed in $\a_{i}$
and for all $i<\k$, $cf(\a_{i+1})= \mu$ and $\a_{i+1}$
is $g$-saturated;

(ii) there are sequence $(\b_{i})_{i<\k}$ and
$\n =\n_{\a}\in J_{0}-J^{-}_{0}$
satisfying (i)-(iv) in the definition of $(J_{0}-J^{-}_{0})
\cup(J_{1}-J^{-}_{1})$.

Let $g(\n)=t(\xi_{0},...,\xi_{n})$, $\xi_{0},...,\xi_{n}\in
J_{1}$.
Now
for all $k\le n$,
either 
$\xi_{k}\in J_{1}(\b_{i})$ for some $i<\k$ or
there is $j<\k$ such that $s(\xi_{k},j)\ge\a$ or
$length(\xi_{k})=\k$, $sup\{ s(\xi_{k},j)\vert\ j<\k\} =\a$ and
for all $j<\k$, $s(\xi_{k},j)<\a$.
By Lemma 3.16, in the last case
$\xi_{k}$ has been put to $J_{1}$ at stage $\a$.

We choose $i<\k$ so that

(a) $i$ is of type 2 and $>2$;

(b) for all $k<l\le n$, $\xi_{k}\raj i\ne\xi_{l}\raj i$;

(c) for all $k\le n$, if
$length(\xi_{k})=\k$, $sup\{ s(\xi_{k},j)\vert\ j<\k\} =\a$ and
for all $j<\k$, $s(\xi_{k},j)<\a$
then there are $\r_{0},...,\r_{r}\in J_{0}\cup
J_{1}$ such that

(i) $\r_{o}=\n$ and $\r_{r}=\xi_{k}$;

(ii) if $p<r$ then $\r_{p}R\r_{p+1}$;

(iii) if $p<r$ then
$W^{\k}_{\r_{p},\r_{p+1}}\subseteq i$;

(iv) for all $p<q\le r$, $\r_{p}\raj i\ne\r_{q}\raj i$;

(d) for all $k\le n$, if $\xi_{k}\in J_{1}(\b_{j})$ for some $j<\k$ then
$\xi_{k}\in J_{1}(\b_{i})$;

(e) for all $k\le n$, if $s(\xi_{k}, j)\ge\a$ for some
$j<\k$ then $\xi_{k}\raj j_{k}\in J_{1}(\b_{i})$ and $j_{k}<i$,
where
$j_{k}=min\{ j<i\vert\ s(\xi_{k},j)\ge\a\}$.

\noindent
Let $l\le l'\le n+1$ be such that
$\xi_{k}\in J_{1}(\b_{i})$ iff $k<l$,
$length(\xi_{k})=\k$, $sup\{ s(\xi_{k},j)\vert$ $j<\k\} =\a$ and
for all $j<\k$, $s(\xi_{k},j)<\a$ iff $l\le k<l'$ and
$\xi_{k}\raj i\not\in J_{1}(\a )$ iff $l'\le k\le n$.
(Of course we may assume that
we have ordered $\xi_{0},...,\xi_{m}$ so that $l$ and $l'$ exist.)
If $l\le k<l'$ then there are
$\r_{0},...,\r_{r}\in J_{1}\cup
J_{0}$ satisfying (c)(i)-(c)(iv) above.
By the choice of $\n (i-1)$,
$\r_{p}\raj i\leftarrow\r_{p+1}\raj i$, for all $p<r$, and
so $\xi_{k}\raj (i+1)\in J_{1}(\b_{i})$.
For all $k\le n$ we define $\xi '_{k}$ as follows:

($\a$) if $k<l$ then $\xi '_{k}=\xi_{k}$;

($\b$) if $l\le k<l'$ then $\xi '_{k}=\xi_{k}\raj (i+1)$;

($\g$) if $l'\le k\le n$ then $\xi '_{k}=\xi_{k}\raj j_{k}$.

Let $g(\n\raj (i+1))=u(\d_{0},...,\d_{m})$, $u$ a term and
$\d_{0},...,\d_{m}\in J_{1}(\b_{i+1})$.
Because $\b_{i}$ is $g$-saturated there is $\n '\in J_{0}(\b_{i})$
and $\d '_{0},...,\d '_{m}\in J_{1}(\b_{i})$ such that

(a) $g(\n ')=u(\d '_{0},...,\d '_{m})$;

(b) $length(\n ')=i+1$ and $\n '\raj i=\n\raj i$;

(c) the basic type of $(\xi '_{0},...,\xi '_{n},
\d_{0},...,\d_{m})$ in
$(J_{1},<,\ll ,H,P_{j})$ is the same as the basic type of
$(\xi '_{0},...,\xi '_{n},
\d '_{0},...,\d '_{m})$.

Because for all $l\le k<l'$, $s(\xi_{k},i+1)\ge\b_{i+1}$ and
for all $l'\le k\le n$, $s(\xi_{k},j_{k})>\b_{i+1}$,
it is easy to see
that the basic type of $(\xi_{0},...,\xi_{n},\d_{0},...,\d_{m})$ in
$(J_{1},<,\ll ,H,P_{j})$ is the same as the basic type of
$(\xi_{0},...,\xi_{n},
\d '_{0},...,\d '_{m})$.

Let $\phi_{n}$, $n<\k$, be as in
Theorem 2.7.
Then
$$EM^{1}(J_{1},\Phi )\models
\phi_{i+1}(u(\d '_{0},...,\d '_{m}),t(\xi_{0},...,\xi_{n})).$$
So $\n '\ne\n\raj (i+1)$, $\n '\raj i=\n\raj i$
and
$$EM^{1}(J_{0},\Phi )\models\phi_{i+1}(\n ',\n ).$$
This is impossible
by Theorem 2.7 (ii). $\eop$

\th 3.19 Conclusion. Let $\l =\mu^{+}$,
$cf(\mu )=\mu$, $\k =cf(\k )<\mu$, $\l^{<\k}=\l$ and $\mu^{\k}=\mu$.
Assume $T$ is an unsuperstable
theory, $\vert T\vert\le\l$ and $\k (T)>\k$. Then there are models
$\A$, $\B\models T$ of cardinality $\l$ such that
$$\A\equiv^{\l}_{\mu\times\k}\B\ \ \hbox{\sl and}
\ \ \A\not\cong\B .$$

\chapter{References}

\item{[HS]} T. Hyttinen and S. Shelah, Constructing strongly equivalent
nonisomorphic models for unsuperstable theories, part A, to appear.

\item{[HT]} T.Hyttinen and H.Tuuri, Constructing strongly equivalent
nonisomorphic models for unstable theories, APAL 52, 1991, 203--248.

\item{[Hu]} T. Huuskonen, Comparing notions of similarity for
uncountable models, Dissertation, University of Helsinki, 1991.

\item{[Sh1]} S.Shelah, Classification Theory, Stud. Logic Found. Math.
92 (North-Holland, Amsterdam, 2nd rev. ed., 1990).

\item{[Sh2]} S.Shelah, Non-structure Theory, to appear.

\bigskip

Tapani Hyttinen

Department of Mathematics

P. O. Box 4

00014 University of Helsinki

Finland
\bigskip

Saharon Shelah

Institute of Mathematics

The Hebrew University

Jerusalem

Israel
\medskip

Rutgers University

Hill Ctr-Busch

New Brunswick

NJ 08903

USA

\end